\documentclass[12pt,oneside,english]{amsart}
\usepackage[latin1]{inputenc}
\usepackage{amssymb}

\makeatletter
\newtheorem{theorem}{Theorem}

\newtheorem{remark}{Remark}
\newtheorem{conjecture}{Conjecture}
\newtheorem{problem}{Problem}

\numberwithin{equation}{section}
\usepackage{babel}

\makeatother
\begin{document}
\baselineskip=17pt

\title[Intervals with at least 3 primes ]{Very small intervals containing at least three primes}

\author{Vladimir Shevelev}
\address{Departments of Mathematics \\Ben-Gurion University of the
 Negev\\Beer-Sheva 84105, Israel. e-mail:shevelev@bgu.ac.il}

\subjclass{Primary 11A41, secondary 11B05}

\begin{abstract}
Let $p_n$ is the $n$-th prime. With help of the Cram\'{e}r-like model, we prove that the set of intervals of the form $(2p_n,\enskip2p_{n+1})$ containing at list 3 primes has a positive density with respect to the set of all intervals of such form.\end{abstract}

\maketitle

\section{Introduction }
 Everywhere below we understand that $p_n$ is the n-th prime and $\mathbb{P}$ is the class of all increasing infinite sequences of primes. If $A\in\mathbb{P}$ then we denote $\mathcal{A}$  the event that prime $p$ is in $A.$ In particular, an important role in our constructions play the following sequences from $\mathbb{P}:\enskip A_i$ is the sequence of those primes $p_k,$ for which the interval $ (2p_k, 2p_{k+1})$ contains at least $i$ primes, $i=1,2,...$  By $\mathcal{A}_i(n),\enskip$ we denote the event that $p_n$ is in $A_i,\enskip i=1,2,...$\newline

In [1] we considered the following problem. Let $p$ be an odd prime. Let, furthermore, $p_n<p/2<p_{n+1}.$ According to the Bertrand's postulate, between $p/2$ and $p$ there exists a prime. Therefore, $p_{n+1}\leq p.$ Again, by the Bertrand's postulate, between $p$ and $2p$ there exists a prime. More subtle question is the following.
 \begin{problem}\label{1} Consider the sequence $S$ of primes $p $ possessing the property: if $p/2$ lies in the interval $(p_n,\enskip p_{n+1})$ then there exists a prime in the interval $(p, 2p_{n+1})$. With what probability a random prime $q$ belongs to $S$ (or the event $\mathcal{S}$ does occur)?.
 \end{problem}
 In this paper we prove the following theorem.
 \begin{theorem}

  The set of intervals of the form $(2p_n,\enskip2p_{n+1})$ containing at list 3 primes has a positive density with respect to the set of all intervals of such form.
   \end{theorem}
   \section{Criterions for $\mathbf{R}$-primes, $\mathbf{L}$-primes and $\mathbf{RL}$-primes }
In [1] we found a sieve for the separating $\mathbf{R}$-primes from all primes and shown how to receive the corresponding sieve for  $\mathbf{L}$-primes. Now we give simple criterions for them.
\begin{theorem}\label{4}
1) $p_n$ is $\mathbf{R}$-prime if and only if $\pi(\frac {p_n} {2})=\pi(\frac {p_{n+1}} {2});$\newline
2) $p_n$ is $\mathbf{L}$-prime if and only if $\pi(\frac {p_n} {2})=\pi(\frac {p_{n-1}} {2})$;\newline
3) $p_n$ is $\mathbf{RL}$-prime if and only if $\pi(\frac {p_{n-1}} {2})=\pi(\frac{p_{n+1}} {2}).$
   \end{theorem}
   \bfseries Proof. \mdseries 1) Let $\pi(\frac {p_n} {2})=\pi(\frac {p_{n+1}} {2})$ is valid. Now if $p_k< p_n/2<p_{k+1},$ and between $p_n/2$ and $p_{n+1}/2$ do not exist primes. Thus $p_{n+1}/2<p_{k+1}$ as well. Therefore, we have $2p_k< p_n<p_{n+1}<2p_{k+1},$ i.e. $p_n$ is $\mathbf{R}$-prime. Conversely, if $p_n$ is $\mathbf{R}$-prime, then $2p_k< p_n<p_{n+1}<2p_{k+1},$ and $\pi(\frac{p_n} {2})=\pi(\frac{p_{n+1}} {2})$ is valid. 2) is proved quite analogously and 3) follows from 1) and 2). $\blacksquare$
  \section{Proof of a "precise symmetry" conjecture   }
   We start with a proof of the following conjecture [1].
    \begin{conjecture}
  Let $\mathbf{R}_n\enskip (\mathbf{L}_n)$ denote the $n$-th term of the sequence $\mathbf{R}\enskip(\mathbf{L}).$
  Then we have
  \begin{equation}\label{2.1}
 \mathbf{R}_1\leq\mathbf{L}_1\leq\mathbf{R}_2\leq\mathbf{L}_2\leq ...\leq\mathbf{R}_n\leq\mathbf{L}_n\leq...
\end{equation}
 \end{conjecture}
 \bfseries Proof of Conjecture 1. \mdseries It is clear that the intervals of considered form, containing not more than one prime, contain neither $\mathbf{R}$-primes nor $\mathbf{L}$-primes. Moving such intervals, consider the first from the remaining ones. The first its prime is an $\mathbf{R}$-prime $(\mathbf{R}_1).$ If it has only two primes, then the second prime is an $\mathbf{L}$-prime $(\mathbf{L}_1), $ and we see that $(\mathbf{R}_1)<(\mathbf{L}_1); $ on the other hand if it has $k$ primes, then beginning with the second one and up to the $(k-1)$-th we have $\mathbf{RL}$-primes, i.e. primes which are simultaneously $\mathbf{R}$-primes and $\mathbf{L}$-primes. Thus, taking into account that the last prime is only $\mathbf{L}$-prime , we have
 $$\mathbf{R}_1<\mathbf{L}_1=\mathbf{R}_2=\mathbf{L}_2=\mathbf{R}_3=...=\mathbf{L}_{k-1}=\mathbf{R}_{k-1}<\mathbf{L}_{k}.$$
 The second remaining interval begins with an  $\mathbf{R}$-prime and the process repeats. $\blacksquare$
 \begin{remark}
 Note that a corollary that "the number of $\mathbf{RL}$-primes not exceeding $x$ is not less than the number of $A_3$-primes not exceeding $x$" is absolutely erroneously. Indeed, we should take into account that every interval of the form $(2p_n,\enskip 2p_{n+1})$ containing $\mathbf{RL}$-prime contains \upshape at least\slshape \enskip3 primes not exceeding $x.$ A right corollary is the following. Since, by the condition of Problem 1, a prime $p$ already lies in a interval $(2p_n,\enskip 2p_{n+1}),$ then we should consider \upshape only \slshape intervals containing at least prime. Denote $\mathcal{A}_{k}, \enskip k=1,...,$ the event that a random interval $(2p_n, \enskip 2p_{n+1})$ contains \slshape at least \upshape $ k,\enskip 1,2,...$ primes. If $P(\mathcal{A}_{1})=q,$ then we have\newpage
 \begin{equation}
 P(\mathcal{A}_{k})=q^k, \enskip k=1,2,...
 \end{equation}
Let, furthermore, $\mathcal{A}^{(k)}, \enskip k=1,...,$ the event that a random interval $(2p_n, \enskip 2p_{n+1})$ contains \slshape exact \upshape $ k,\enskip 1,2,...$ primes. Then, by (3.2),
$$P(\mathcal{A}^{(k)}= P(\mathcal{A}_{k})- P(\mathcal{A}_{k+1})=(q-1)q^k, \enskip k=1,2,...$$
and we have
 \begin{equation}
 P(\mathbf{RL})=(1-q)\sum_{k\geq3}\frac {k-2} {k}q^{k-1}=2-q+2\frac{1-q} {q}\ln(1-q).
 \end{equation}
\end{remark}

\section{Proof of Theorem 1}
The theorem immediately follows from the positivity of probability $P(\mathbf{RL}).$ In fact, in [1] we proved that $q\approx0.8010$ and $P(\mathbf{RL})\approx0.3980.$
 $\blacksquare$\newline
Note that by the Cram\'{e}r's 1937 conjecture $2p_{n+1}-2p_n<(2+\varepsilon)\ln^2 n.$ Thus, there exists an infinite sequence of the intervals of such small length, but having at least three primes, and, moreover, this sequence has
a positive density with respect to the sequence of all intervals of the form $(2p_n,\enskip 2p_{n+1}).$\newline
By this way, in view of (3.2), it could be proved a more general result.
\begin{theorem}

  Let $h$ be arbitrary large but a fixed positive integer. Then the set of intervals of the form $(2p_n,\enskip2p_{n+1})$ containing at list $h$ primes has a positive density with respect to the set of all intervals of such form.
   \end{theorem}

Quite analogously one can consider an $m$-generalization of Theorem 1 for every $1<m<2.$ Here the case of especial interest is the case of the values of $m$ close to 1.

\;\;\;\;\;\;\;\;


\begin{thebibliography}{1}
\bibitem 1. V. Shevelev \slshape Three probabilities concerning prime gaps \upshape http://arxiv.org/abs/0909.0715
\bibitem 2. V. Shevelev \slshape Critical small intervals containing primes\upshape http://arxiv.org/abs/0908.2319
\end{thebibliography}
\end{document}